\documentclass{ifacconf}

\usepackage{graphicx}      
\usepackage{natbib}        

\def\be{\begin{equation}}
\def\ee{\end{equation}}
\def\ba{\begin{array}}
\def\ea{\end{array}}
\usepackage[gen]{eurosym}
\usepackage{epstopdf}
\usepackage{dsfont}

\usepackage{multicol}

\usepackage{pgf}
\usepackage{tikz}
\usetikzlibrary{arrows,automata}
\usepackage{amsmath}
\usetikzlibrary{arrows,shadows,calc,decorations.shapes,decorations} \usepackage[underline=true,rounded corners=false]{pgf-umlsd}
\usepackage[american,cute inductors,smartlabels]{circuitikz}
\usepackage{mathrsfs}
\usepackage{subfigure}
\usepackage{graphicx,color}
\usepackage{graphics} 	
\usepackage{latexsym}
\usepackage{booktabs}
\usepackage{amssymb}
\usepackage{amsmath,amssymb,amsfonts}

\DeclareMathOperator{\diag}{diag}
\usepackage{mathtools}

\newtheorem{remark}{Remark}
\newtheorem{goal}{Goal}
\newtheorem{definition}{Definition}

\newtheorem{theorem}{Theorem}

\newtheorem{lemma}{Lemma}

\newtheorem{assumption}{Assumption}

\usepackage{siunitx}	
\usepackage{balance}										 		
\DeclareMathAlphabet{\mathpzc}{OT1}{pzc}{m}{it}
\DeclareMathAlphabet{\mathcal}{OMS}{cmsy}{m}{n}

\newcommand{\R}{\mathbb{R}}

\begin{document}

\begin{frontmatter}

\title{Passivity properties for regulation of DC networks with stochastic load demand
}
\author[First,Second]{Amirreza Silani} 
\author[First]{Michele Cucuzzella} 
\author[First]{Jacquelien M. A. Scherpen}
\author[Second]{Mohammad Javad Yazdanpanah}

\address[First]{Jan C. Wilems Center for Systems and Control, ENTEG, Faculty of Science and Engineering, University of Groningen, Nijenborgh 4, 9747 AG Groningen, the Netherlands. (e-mail: \{a.silani, m.cucuzzella, j.m.a.scherpen\}@rug.nl).}
\address[Second]{Control \& Intelligent Processing Center of Excellence, School of Electrical and Computer Engineering, University of Tehran, Tehran, Iran (e-mail: yazdan@ut.ac.ir)}

%

\begin{abstract}                
In this paper we present new (stochastic) passivity properties for Direct Current (DC) power networks, where the unknown and unpredictable load demand is modelled by a stochastic process. 
More precisely, the considered power network consists of distributed generation units supplying ZIP loads, \textit{i.e.}, nonlinear loads comprised of impedance (Z), current (I) and power (P) components. Differently from the majority of the results in the literature, where each of these components is assumed to be constant, we consider time-varying loads whose dynamics are described by a class of stochastic differential equations. 
Finally, we prove that an existing distributed control scheme achieving current sharing and (average) voltage regulation ensures the asymptotic stochastic stability of the controlled network.
\end{abstract}


\end{frontmatter}
\section{INTRODUCTION}
Power networks are typically classified into Direct Current (DC) and Alternating Current (AC) networks, which are interconnected clusters of Distributed Generation Units (DGUs), loads and energy storage devices. The recent wide spread of renewable energy sources, electronic appliances and batteries (including electric vehicles) motivates the design and operation of DC networks, which are generally more effcient and reliable than AC networks (\cite{ref1}), attracting growing interest.

In order to guarantee a proper and safe functioning of the overall network and the appliances connected to it, the main goal in DC networks is voltage stabilization (see for instance \cite{ref2,ref17,ref3,ref23,Iovine,ref13,ref22,ref24,CucuzzellaCDC19}). Moreover, as different DGUs may generally have different generation (or storage) capacities, an additional goal is to (fairly) share the total demand of the network among its DGUs (see for instance \cite{ref8,ref12,ref24.5,ref20.5,Martinelli}). This goal is usually called power or current sharing and its achievement does not generally permit to regulate the voltage at each node towards the corresponding pre-specified reference value. Consequently, different forms of voltage regulation have been proposed in the literature, where for instance the average value of the voltages of the whole microgrid is controlled towards a desired setpoint (see for instance \cite{ref8,ref12,ref24.5}).


\subsection{Motivation and contributions}
In all these works the load components are assumed to be constant. 
However, it is well known that electric loads are in practice time-varying and, due to the random and unpredictable diversity of usage patterns, it is more realistic to consider unknown time-varying loads described for instance by stochastic processes (see for instance \cite{ref24.55,ref24.10,ref31}). In \cite{ref24.6}, a control scheme is proposed to regulate the renewable energy sources and energy storage devices of a DC network with stochastic Z loads. The control strategy requires the design of a filter for estimating the load resistance.
A cascade control system for the energy management of DC microgrids with I loads is presented in \cite{ref24.8}, where the proposed control scheme includes an adaptive estimation of the quasi-stochastic load current profiles. 
In \cite{ref24.9}, a droop control scheme is designed for DC microgrids with stochastic Z loads. {Moreover, in some papers, the Stochastic Differential Equations (SDEs) have been used for modeling the loads and other uncertainties in power system networks (see for instance \cite{new4,new6,new7}). In \cite{new4}, the random load characteristic is considered to develop a stochastic model for voltage stability analysis. A stochastic power system model based on stochastic differential equations is presented in \cite{new6} to consider the uncertain factors such as load levels and system faults. In \cite{new7}, a systematic and general
approach to model power systems as continuous stochastic differential-algebraic equations is proposed and it justifies the need for stochastic models in power system analysis.}

In this paper, (i) each component of the ZIP load is modeled as the sum of an unknown constant and the solution to a stochastic differential equation  describing the load dynamics; (ii) sufficient conditions for the stochastic passivity of the open-loop system are presented, facilitating the interconnection with passive control systems; (iii) the asymptotic stochastic stability of the power network controlled by the distributed control scheme proposed by \cite{ref24.5} is proved. 

%
%
%

\subsection{Outline}

The paper is organized as follows. The DC network model is introduced in Section~2. In Section~3, the control problem is formulated. In Section~4, the stochastic dynamics of the ZIP loads are introduced and the stochastic passivity properties of the considered DC network are shown, then the stochastic stability of the closed-loop system with an existing controller is proved. The simulation results are presented and discussed in Section~5, while some conclusions and future research directions are gathered in Section~6. 

\subsection{Notation}
The set of real numbers is denoted by $\R$. The set of positive (nonnegative) real numbers is denoted by $\R_{>0}$ ($\R_{\geq0}$).
Let $\vec{0}$ be the vector of all zeros or the null matrix of suitable dimension(s) and let $\mathds{1}_n \in \R^n$ be the vector containing all ones. The $i$-th element of vector $x$ is denoted by $x_i$. 
 A steady-state solution to the system $\dot x = \zeta(x)$, is denoted by $\overline x$,  \textit{i.e.}, $\boldsymbol{0} = \zeta(\overline x)$. 
Moreover, given a signal (or parameter) $\sigma$, $\sigma^\ast$ denotes that $\sigma$ is a constant signal (or parameter), while $\hat\sigma$ denotes that $\sigma$ is a stochastic process. Given a vector $x\in \R^n$, $[x]\in \R^{n\times n}$ indicates the diagonal matrix whose diagonal entries are the components of $x$. 
 Let $A \in \R^{n \times n}$ be a matrix. In case $A$ is a positive definite (positive semi-definite) matrix, we write $A > \vec{0}$  ($A \geq \vec{0}$). The $n \times n$ identity matrix is denoted by  $\mathbb{I}_{n}$ and the vector consisting of all ones is denoted by $\mathds{1}_{n}$.
 
\section{DC Network Model}
\begin{figure}[t]
    	\begin{center}
	\ctikzset{bipoles/length=1cm}
    		\begin{circuitikz}[scale=1,transform shape]
    			\ctikzset{current/distance=1}
    			\draw
    			node[] (Ti) at (0,0) {}
    			node[] (Tj) at ($(5.4,0)$) {}
    			node[] (Aibattery) at ([xshift=-4.5cm,yshift=0.9cm]Ti) {}
    			node[] (Bibattery) at ([xshift=-4.5cm,yshift=-0.9cm]Ti) {}
    			node[] (Ai) at ($(Aibattery)+(0,0.2)$) {}
    			node[] (Bi) at ($(Bibattery)+(0,-0.2)$) {}
    			($(Ai)$) to [short,i={$I_{gi}$}]($(Ai)+(.2,0)$){}
    			($(Ai)+(.2,0)$) to [L, l={$L_{gi}$}] ($(Ai)+(1.5,0)$){}
    			to [short, l={}]($(Ti)+(0,1.1)$){}
    			(Bi) to [short] ($(Ti)+(0,-1.1)$);
    			\draw
    			($(Ai)$) to []($(Aibattery)+(0,0)$)to [V_=$u_i$]($(Bi)$)
    			($(Ti)+(-2.8,1.1)$) node[anchor=south]{{$V_{i}$}}
    			($(Ti)+(-2.8,1.1)$) node[ocirc](PCCi){}
    			($(Ti)+(-1.8,1.1)$) to [R, l={$G_{li}$}] ($(Ti)+(-1.8,-1.1)$)
    			($(Ti)+(-.6,1.1)$) to [short,i>={$I_{li}$}]($(Ti)+(-.6,0.5)$)to [I]($(Ti)+(-.6,-1.1)$)
    			($(Ti)+(.2,1.1)$) to [short,i>={$\dfrac{P_{li}}{V_i}$}]($(Ti)+(.2,0.5)$)to [I]($(Ti)+(.2,-1.1)$)
    			($(Ti)+(-2.8,1.1)$) to [C, l_={$C_{gi}$}] ($(Ti)+(-2.8,-1.1)$)
    			($(Ti)+(2.,1.1)$) to [short,i={$I_{k}$}] ($(Ti)+(2.2,1.1)$)
    			($(Ti)+(0,1.1)$)--($(Ti)+(.6,1.1)$) to [R, l={$R_{k}$}] 
    			($(Ti)+(2.5,1.1)$) {} to [L, l={{$L_{k}$}}, color=black]($(Tj)+(-2.2,1.1)$){}
    			($(Tj)+(-2.2,1.1)$) to [short]  ($(Ti)+(3.4,1.1)$)
    			($(Ti)+(0,-1.1)$) to [short] ($(Ti)+(3.4,-1.1)$);
    			\draw
    			node [rectangle,draw,minimum width=2.9cm,minimum height=3.4cm,dashed,color=gray,label=\textbf{DGU $i$},densely dashed, rounded corners] (DGUi) at ($0.5*(Aibattery)+0.5*(Bibattery)+(.65,0.2)$) {}
    			node [rectangle,draw,minimum width=3.1cm,minimum height=3.4cm,dashed,color=gray,label=\textbf{ZIP $i$},densely dashed, rounded corners] (DGUi) at ($0.5*(Aibattery)+0.5*(Bibattery)+(3.825,0.2)$) {}
    			node [rectangle,draw,minimum width=2.2cm,minimum height=3.4cm,dashed,color=gray,label=\textbf{Line $k$},densely dashed, rounded corners] (DGUi) at ($0.5*(Aibattery)+0.5*(Bibattery)+(6.65,0.2)$) {};
    		\end{circuitikz}
    		\caption{Electrical scheme of DGU $i$, ZIP load $i$ and transmission line $k$, with $i\in\mathcal{V}$ and $k\in\mathcal{E}$.}
    		\label{f1}
    	\end{center}
    \end{figure}
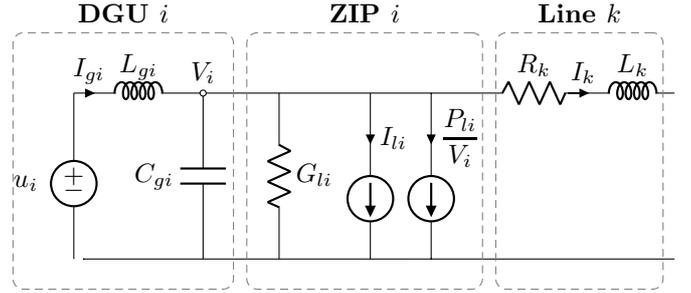
    
    \begin{table}
\centering
\caption{Description of symbols}
{\begin{tabular}{ll}
\midrule
\midrule
${G}_{li}$		& Unknown conductance load\\
${I}_{li}$		& Unknown current load\\  
${P}_{li}$		& Unknown power load\\  
$I_{gi}$						& Generated current\\
$V_i$						    & Load voltage\\
$I_{k}$ 						& Line current \\
$u_i$						& Control input \\
${L_{k}}$	     		& {Line inductance}\\
$R_{k}$						& Line resistance\\
$C_{gi}$						& Shunt capacitor\\
$L_{gi}$						& Filter inductance\\
\midrule
\midrule
\end{tabular}}
\label{tab1}
\end{table}

In this section, we introduce the model of the considered DC network comprising of Distributed Generation Units (DGUs), loads and transmission lines (\cite{ref8,ref12,ref24.5}). Fig.~{\ref{f1}} illustrates the structure of the considered DC microgrid and Table~\ref{tab1} reports the description of the used symbols. 
Let $\mathcal{G} = (\mathcal{V},\mathcal{E})$ be a connected and undirected graph that describes the considered DC network. The nodes and the edges are denoted by $\mathcal{V} = \{1,...,n\}$ and $\mathcal{E}  = \{1,...,m\}$, respectively. Then, the  topology of the microgrid is represented by the corresponding incidence matrix $\mathcal{A} \in \R^{n \times m}$. Let the ends of each transmission line $k$ be arbitrarily labeled with a `-', or a `+', then $\mathcal{A}$ can be defined as 
\begin{equation}\label{eq3}
\mathcal{A}_{ik}=
\begin{cases}
+1 \quad &\text{if $i$ is the positive end of $k$}\\
-1 \quad &\text{if $i$ is the negative end of $k$}\\
0 \quad &\text{otherwise}.
\end{cases}
\end{equation}
Before presenting the dynamics of the overall network, we first introduce and discuss the models of each component of the network.

{\bf DGU model:} The dynamic model of DGU $i\in\mathcal{V}$ is described by
\begin{align}
\begin{split}
\label{eq1}
L_{gi}\dot{I}_{gi} &=  -  V_i + u_i \\
C_{gi}\dot{V}_i &=   I_{gi} - {I}_{li}(V_i) - \displaystyle{ \sum_{k \in \mathcal{E}_i}^{}I_{k}},
\end{split}
\end{align}
where $I_{gi}, V_i, I_{k}, u_i: \R_{\geq 0} \rightarrow\R$, $L_{gi}, C_{gi} \in \R_{>0}$ and $\mathcal{E}_i$ is the set of the lines incident to node $i$. Moreover, $I_{li}: \R \rightarrow \R$ represents the current demand of load $i$ possibly depending on the voltage $V_i$.

{\bf Load model:} In this work, we consider a general \emph{nonlinear} and \emph{stochastic} load model including the parallel combination of the following \emph{unknown} load components\footnote{Note that for the sake of notation simplicity we prefer to use the load conductance $G_{li}$ instead of the load impedance $Z_{li}=G_{li}^{-1}$.}:
\begin{enumerate}
\item{impedance component (Z): $G^\ast_{li} + \hat{G}_{li}$,}
\item{current component (I): $I^\ast_{li}+\hat{I}_{li}$,}
\item{power component (P): $P^\ast_{li}+\hat{P}_{li}$,}
\end{enumerate}
where $G^\ast_{li}, I^\ast_{li}, P^\ast_{li} \in\R_{>0}$ are \emph{unknown} constants and $\hat{G}_{li}, \hat{I}_{li},$ $ \hat{P}_{li}: \R_{\geq 0} \rightarrow\R$ are stochastic processes, the dynamics of which will be introduced in Section \ref{sec:passivity}.
To refer to the load types above, the letters Z, I and P, respectively, are often used in the literature (see for instance \cite{ref20.5}).
 We also use Z$^\ast$, I$^\ast$ and P$^\ast$ to denote the unknown constant components of the load\footnote{For instance, in presence of a Z$^\ast$IP$^\ast$ load, we have $I_{li}(V_i) = G^\ast_{li}V_i + I^\ast_{li}+\hat{I}_{li} + V_i^{-1}P^\ast_{li}$.}. 
Therefore, in presence of ZIP loads, $I_{li}(V_i)$ in \eqref{eq1} is given by
\begin{equation}
\label{eq:plant_i3}
I_{li}(V_i) = (G^\ast _{li} + \hat{G}_{li})V_i + I^\ast_{li}+\hat{I}_{li} + V_i^{-1}(P^\ast_{li}+\hat{P}_{li}).
\end{equation}

{\bf Line model:} The dynamics of the current $I_{k}$ exchanged between nodes $i$ and $j$ are described by 
\begin{equation}
{
\label{eq2}
L_{k}{\dot I_{k}} = (V_i - V_j) - R_{k} I_{k},}
\end{equation}
where  $R_k, L_k \in \R_{>0}$.

Now, the dynamics of the overall network can be written compactly as 
\begin{equation}
\label{eq4}
	\begin{split}
		L_g\dot{I}_{g} & =- V + u\\
		C_g\dot{V} & = I_{g} + \mathcal{A}I - [G_{l}]V - I_{l} - [V]^{-1}P_{l}\\
		L\dot{I} & = -\mathcal{A}^\top V - R I,
	\end{split}
\end{equation}
where $I_g, V, u : \R_{\geq 0}\rightarrow\R^{n}$, $I: \R_{\geq 0}\rightarrow\R^m$, $ L_g, C_g \in \R^{n \times n}_{>0}$ and $R, L \in \R^{m \times m}_{>0}$. Moreover, $G_l = G^\ast_{l} + \hat{G}_{l}, I_l=I^\ast_{l}+\hat{I}_{l}$ and $P_l = P^\ast_{l}+\hat{P}_{l}$, where $G_l^\ast, I_l^\ast, P_l^\ast \in \R^n_{>0}$ represent the constant components of the load and $\hat{G}_l,  \hat{I}_l, \hat{P}_l:\R_{\geq 0} \rightarrow \R^n$ the stochastic parts. These dynamics are discussed in Section 4, where we assume that their steady-state solutions are equal to zero.

\section{Problem Formulation}
In this section, we introduce the main control goals of this paper, \textit{i.e.}, current sharing and average voltage regulation. Assuming that the stochastic components of the ZIP loads are zero at the steady-state, we notice that for a constant input $u=u^\ast$, the steady-state solution $(\bar{I}_g,\bar{V},\bar{I})$ to \eqref{eq4} satisfies
\begin{subequations}\label{eq5}
\begin{align}
		\overline{V} & = u^\ast \label{ss1}\\
		 I_l^\ast+[G_l^\ast]\bar{V}+[\bar{V}]^{-1}P_l^\ast -  \overline{I}_{g} &= \mathcal{A} \overline I \label{ss2}\\
		\overline I & = -R^{-1}\mathcal{A}^\top\overline V. \label{ss3}
\end{align}
\end{subequations}

Because of different generation capacities, it is reasonable to require that the total load demand of the microgrid is fairly shared among all the different generation units as $q_i \bar{I}_{gi}=q_j\bar{I}_{gj},~ \forall i,j\in \mathcal{V}$, where $q_i$ has inverse relationship with the capacity of DGU $i$. Then, we define the first goal concerning the steady-state value of the generated current as follows:
\begin{goal}{\bf(Current sharing).}
\label{obj:current_sharing}
 \begin{align}
 \label{eq6}
\lim_{t \rightarrow \infty} I_g(t)= \overline{I}_{g} = Q^{-1}\mathds{1}_{n} {i}_{g}^\ast,
\end{align}
where $Q= \diag(q_1,\dots,q_n)$, $q_i \in \R_{>0},~ \forall i \in \mathcal{V}$ and $i_g^\ast$ any constant value satisfying at the steady-state ${i}_{g}^\ast = \frac{\mathds{1}_{n}^\top ([G_l^\ast] \bar{V}+{I}_l^\ast+[\bar{V}]^{-1}P_l^\ast)}{\mathds{1}_{n}^\top Q^{-1}\mathds{1}_{n}}$. Indeed, the latter implies that at the steady-state the total generated current is equal to the total load demand, \textit{i.e.}, $\mathds{1}_{n}^\top\bar{I}_g=\mathds{1}_{n}^\top  I_l(\bar{V})$.
\end{goal}

\smallskip

Now, we observe that achieving Goal 1, does not generally allow to perform also voltage regulation. Indeed, the achievement of current sharing generally implies voltage differences between the nodes of the  microgrid. As a consequence, a particular form of voltage regulation has been proposed in the literature, where the (weighted) average of the steady-state voltages is regulated towards the (weighted) average of the voltage references (\cite{ref12,ref24.5}). Then, the second goal concerning the steady-state value of the voltages is defined as follows: 
\smallskip
\begin{goal}{\bf (Average voltage regulation).}
\label{obj:voltage_balancing}
\begin{align}
\label{eq:voltage_balancing}
\lim_{t \rightarrow \infty} {\mathds{1}_{n}^\top Q^{-1} V(t)} =  {\mathds{1}_{n}^\top Q^{-1} \overline{V}} = {\mathds{1}_{n}^\top Q^{-1}V^{\ast}},
\end{align}
$V^{\ast}_i$ being the voltage reference at node $i\in\mathcal{V}$.
\end{goal}

Before showing the stochastic properties of the controlled network, we assume that there exists a steady-state solution to \eqref{eq4}:

\begin{assumption}{\bf (Steady-state solution). }\label{ass:1}
There exists a constant input $u^{\ast}$ and a steady-state solution $(\bar{I}_g,\bar{V},\bar{I})$ to \eqref{eq4} satisfying \eqref{eq5} and achieving Goals 1 and 2.
\end{assumption}

\section{Stochastic Passivity of DC Networks}
\label{sec:passivity}
In this section, we introduce the dynamics of the stochastic components of the load, \textit{i.e.}, $\hat{I}_{l}$, $\hat{P}_l$ and $\hat{G}_{l}$, respectively, and verify the stochastic passivity of the open-loop system. Finally, we prove that the control scheme proposed by \cite{ref24.5} ensures the asymptotic stochastic stability of the controlled network.

First, we notice that the stochastic terms in \eqref{eq4} enter in a multiplicative manner.
	For this reason, an appropriate mathematical framework such as the Ito calculus framework, should be adopted
	to analyze such a model. Therefore, we recall for the readers' convenience the definitions of stochastic differential equation and stochastic passivity through the Ito calculus framework (\cite{ref25,ref26,ref26.1,ref30}).

\begin{definition}{\bf (Stochastic differential equation).} We define a stochastic differential equation (SDE) as
\begin{equation}\label{eq7}
dx(t)=f(x,u)dt+g(x)dw(t),
\end{equation}
where $f(x,u)\in\mathbb{R}^{N}$ and $g(x)\in\mathbb{R}^{N\times M}$ are locally Lipschitz, $x(t)\in\mathbb{R}^N$ is the state vector, $u(t)\in\mathbb{R}^P$ is the input of the system and $w(t)\in \mathbb{R}^{M}$ is the standard Brownian motion vector. In the following, we define the Ito derivative (\cite{ref25,ref26,ref26.1,ref30}).
\end{definition}
\begin{definition}{\bf (Ito derivative).} Consider a storage function $S(x)$, which is twice continuously differentiable. Then, $\mathcal{L}S(x)$ denotes the Ito derivative of $S(x)$ along the SDE (\ref{eq7}), \textit{i.e.},
\begin{equation}
\begin{split}\label{eq8}
\mathcal{L}S(x)=\dfrac{\partial S(x)}{\partial x}f(x,u)+\frac{1}{2}\text{tr}\{g^\top (x)\dfrac{\partial ^2 S(x)}{\partial x^\top \partial x}g(x)\}.
\end{split}
\end{equation}
\end{definition}
Now, for the readers' convenience, we introduce the concept of stochastic passivity (\cite{ref26,ref26.1,ref30}).
\begin{definition} {\bf (Stochastic passivity).}~Consider system (\ref{eq7}) with output $y=h(x)$. Assume that the deterministic and stochastic terms of the SDE (\ref{eq7}) at the equilibrium point are identical to zero, \textit{i.e.}, $f(\bar{x},{u^\ast})=g(\bar{x})=\textbf{0}$. Then, system (\ref{eq7}) is said to be stochastically passive with respect to the supply rate $u^\top y$ if there exists a twice continuously differentiable positive semi-definite storage function $S(x)$ satisfying
\begin{equation}\label{eq10}
\mathcal{L}S(x)\leq u^\top y,~ \forall (x,u) \in\mathbb{R}^N\times\mathbb{R}^P.
\end{equation}
\end{definition}

Now, in the following subsections we introduce the dynamics of the stochastic load components, which we model via an SDE in the Ito calculus framework. More precisely, we adopt the well-known Ornstein-Uhlenbeck process (\cite{refss}), which is indeed widely used for the description of physical phenomena (see for instance \cite{ref31}).
For the sake of exposition, we first consider Z$^\ast$IP$^\ast$ loads, \textit{i.e.}, only the load current is described by a stochastic process. Later, we extend the results also to Z$^\ast$IP and ZIP loads, respectively.

\subsection{Z$^\ast$IP$^\ast$ loads}
In this subsection, we consider Z$^\ast$IP$^\ast$ loads type, \textit{i.e.}, equation \eqref{eq:plant_i3} becomes $I_{li}(V_i) = G^\ast_{li}V_i+ I^\ast_{li}+\hat{I}_{li} + V_i^{-1}P^\ast_{li}$, where the SDE describing the dynamics of $\hat{I}_{li}$ is given by 
\begin{equation} \label{eq11}
d\hat{I}_{li}=-\mu _{I_{li}}\hat{I}_{li} dt + \sigma _{I_{li}} \hat{I}_{li} dw_i,
\end{equation}
where $\mu _{I_{li}}$ and $\sigma _{I_{li}}$ are positive constants. As mentioned in Section~3, it can be noticed that the deterministic and stochastic parts of (\ref{eq11}) at the equilibrium point are identical to zero.
Now, in order to verify a passivity property of the considered network, we introduce an assumption on the parameters of the SDE \eqref{eq11} and restrict the trajectories to be inside a subset of the state-space.
\begin{assumption}{\bf(Condition on the parameters of \eqref{eq11}).}\label{ass:2}
For all $i\in\mathcal{V}$, the stochastic parameters of (\ref{eq11}) satisfy
\begin{equation}\label{eq13}
\mu _{I_l} >\frac{1}{2}\sigma _{I_l}^2-\frac{1}{2} \mathbb{I}_{n},
\end{equation}
where $\mu _{I_l}=\diag\{\mu _{I_{l1}}, \dots, \mu _{I_{ln}}\}$ and $\sigma _{I_l}=\diag\{\sigma _{I_{l1}}, \dots,$ $\sigma _{I_{ln}}\}$.
\end{assumption}
Let us define the set $\Omega_\mathrm{Z^\ast IP^\ast}:=\{ (I_g,V,I,\hat{I}_l): \frac{1}{2}([G_l^\ast]-\Pi^{-1})-[P^\ast_l][V]^{-1}[\bar{V}]^{-1} > 0\}$.
Then, the stochastic passivity property of system (\ref{eq4}) with Z$^\ast$IP$^{\ast}$ loads is verified via the following lemma:
\begin{lemma}{\bf (Stochastic passivity of \eqref{eq4}, \eqref{eq11}).}
Let Assumption \ref{ass:2} hold and $\Omega_\mathrm{Z^\ast IP^\ast}$ be nonempty. System (\ref{eq4}) with Z$^{\ast}$IP$^{\ast}$ loads and $\hat{I}_l$ given by (\ref{eq11}) is stochastically (shifted) passive with respect to the storage function 
\begin{equation}\label{eq14}
\begin{split}
S_\mathrm{Z^\ast IP^\ast}(I_g, V, I, \hat{I}_l)=& \frac{1}{2} (I_g -\bar{I}_g)^\top L_g(I_g-\bar{I}_g)\\
&+\frac{1}{2}(V-\bar{V})^\top C_g(V-\bar{V})\\
&+\frac{1}{2}(I-\bar{I})^\top L(I-\bar{I})+\frac{1}{2}\hat{I}_l^\top\Pi\hat{I}_l,
\end{split}
\end{equation}
and supply rate $(I_g-\bar{I}_g)^\top (u-{u^\ast})$ for all the trajectories $(I_g,V,I,\hat{I}_l)\in\Omega_\mathrm{Z^\ast I P^\ast}$, where $\Pi$ is a (suitable) positive definite diagonal constant matrix of appropriate dimensions and $(\bar{I}_g,\bar{V},\bar{I})$ satisfies \eqref{eq5}.
\end{lemma}
\begin{pf}
The Ito derivative of the storage function (\ref{eq14}) satisfies
\begin{equation}\label{eq15}
\begin{split}
\mathcal{L}S_\mathrm{Z^\ast IP^\ast}=
&(I_g-\bar{I}_g)^\top (u-{u^\ast})-(I-\bar{I})^\top R(I-\bar{I})\\
&-\hat{I}_l^\top (\mu_{I_l}-\frac{1}{2} \sigma _{I_l}^2-\frac{1}{2}\mathbb{I}_{n})\Pi\hat{I}_l\\
&-\frac{1}{2}(\Pi ^{-\frac{1}{2}}(V-\bar{V})+\Pi ^{\frac{1}{2}}\hat{I}_l)^\top (\Pi^{-\frac{1}{2}}(V-\bar{V})\\
&+\Pi^{\frac{1}{2}}\hat{I}_l)-(V-\bar{V})^\top (\frac{1}{2}([G^\ast _l]-\Pi ^{-1})\\
&-[P^\ast _l][V]^{-1}[\bar{V}]^{-1})(V-\bar{V}),
\end{split}
\end{equation}
along the solutions to \eqref{eq4}, \eqref{eq11}. Then, we can conclude that $\mathcal{L}S_\mathrm{Z^\ast IP^\ast}(x) \leq  (u-{u^\ast})^\top (I_g-\bar{I}_g)$ for all the trajectories $(I_g,V,I,\hat{I}_l)\in \Omega_\mathrm{Z^\ast IP^\ast}$.
\hfill $\blacksquare$
\end{pf}

\begin{remark}{\bf (Z$^\ast$I load).}
We observe that in presence of only Z$^\ast$I loads, the result provided in Lemma~1 can be strengthened. Indeed, the absence of P loads implies that the system \eqref{eq4}, \eqref{eq11} is stochastically passive for any sufficiently large $\Pi$.
\end{remark}
\subsection{Z$^\ast$IP loads}
In this subsection, we consider Z$^\ast$IP loads type, \textit{i.e.}, equation \eqref{eq:plant_i3} becomes $I_{li}(V_i) = G^\ast_{li}V_i + I^\ast_{li}+\hat{I}_{li} + V_i^{-1}(P^\ast_{li}+\hat{P}_{li})$, where the SDE describing the dynamics of $\hat{P}_{li}$ is given by 
\begin{equation} \label{eq16}
d\hat{P}_{li}=-\mu _{P_{li}}\hat{P}_{li} dt + \sigma _{P_{li}} \hat{P}_{li} dw_i,
\end{equation}
where $\mu _{P_{li}}$ and $\sigma _{P_{li}}$ are positive constants. Also in this case, we notice that the equilibrium points of the deterministic and stochastic parts of (\ref{eq16}) are zero.
Now, in order to verify a passivity property of the considered network, we introduce an assumption on the parameters of the SDE \eqref{eq16} and restrict the trajectories to be inside a subset of the state-space. 
\begin{assumption}{\bf(Condition on the parameters of \eqref{eq16}).}\label{ass:4}
For all $i\in\mathcal{V}$, the stochastic parameters of (\ref{eq16}) satisfy
\begin{equation}\label{eq17}
\mu _{P_l} >\sigma _{P_l}^2,
\end{equation}
\end{assumption}
where $\mu _{P_l}=\diag\{\mu _{P_{l1}},\dots,\mu _{P_{ln}} \}$ and $\sigma _{P_l}=\diag\{\sigma _{P_{l1}},\dots,$ $\sigma _{P_{ln}}\}$. 

Let us define the set $\Omega_\mathrm{Z^\ast IP}:=\{ (I_g,V,I,\hat{I}_l,$ $\hat{P}_l): \frac{1}{2}([G_l^\ast]-\Pi^{-1})-[P^\ast_l][V]^{-1}$ $[\bar{V}]^{-1}-\frac{1}{2}\Sigma\mu _{P_l}[V]^{-2} > 0\}$.
Then, the stochastic passivity property of system (\ref{eq4}) with Z$^\ast$IP loads is verified via the following lemma:
\begin{lemma}{\bf (Stochastic passivity of  \eqref{eq4}, \eqref{eq11}, \eqref{eq16}).}
Let Assumptions \ref{ass:2}, \ref{ass:4} hold and $\Omega_\mathrm{Z^\ast IP}$ be nonempty. System (\ref{eq4}) with Z$^\ast$IP loads and $\hat{I}_l, \hat{P}_l$ given by (\ref{eq11}), (\ref{eq16}) is stochastically (shifted) passive with respect to the storage function 
\begin{equation}\label{eq18.1}
\begin{split}
S_\mathrm{Z^\ast IP}(I_g, V, I, \hat{I}_l,\hat{P}_l)=&~S_\mathrm{Z^\ast IP^\ast}(I_g, V, I, \hat{I}_l)+\frac{1}{2}\hat{P}_l^\top\Sigma \hat{P}_l,
\end{split}
\end{equation}
and supply rate $(I_g-\bar{I}_g)^\top (u-{u^\ast})$ for all the trajectories $(I_g,V,I,\hat{I}_l,\hat{P}_l)\in\Omega_\mathrm{Z^\ast IP}$, where $\Sigma$ is a (suitable) positive-definite diagonal constant matrix of appropriate dimensions.
\end{lemma}
\begin{pf}
The Ito derivative of the storage function (\ref{eq18.1}) satisfies
\begin{equation}\label{eq19}
\begin{split}
\mathcal{L}S_\mathrm{Z^\ast IP}=
&(I_g-\bar{I}_g)^\top (u-{u^\ast})-(I-\bar{I})^\top R(I-\bar{I})-\hat{I}_l^\top (\mu_{I_l}\\
&-\frac{1}{2} \sigma _{I_l}^2-\frac{1}{2}\mathbb{I}_{n})\Pi\hat{I}_l-\frac{1}{2}\big (\Pi^{-\frac{1}{2}}(V-\bar{V})\\
&+\Pi^{\frac{1}{2}}\hat{I}_l\big )^\top\big (\Pi^{-\frac{1}{2}}(V-\bar{V})+\Pi^{\frac{1}{2}}\hat{I}_l\big )\\
&-\frac{1}{2}\hat{P}_l^\top (\mu _{P_l}-\sigma _{P_l}^2)\Sigma \hat{P}_l-\frac{1}{2}\Big (\Sigma ^{\frac{1}{2}}\mu _{P_l}^\frac{1}{2}\hat{P}_l\\
&-\Sigma ^{-\frac{1}{2}}\mu _{P_l}^{-\frac{1}{2}}(-\mathds{1}_{n}^\top +\bar{V}^\top [V]^{-1})^\top\Big )^\top \\
&\Big (\Sigma ^{\frac{1}{2}}\mu _{P_l}^\frac{1}{2}\hat{P}_l-\Sigma ^{-\frac{1}{2}}\mu _{P_l}^{-\frac{1}{2}}(-\mathds{1}_{n}^\top +\bar{V}^\top [V]^{-1})^\top\Big )\\
&-(V-\bar{V})^\top\Big (\frac{1}{2}([G^\ast _l]-\Pi ^{-1})-[P^\ast_l][V]^{-1}[\bar{V}]^{-1}\\
&-\frac{1}{2}\Sigma\mu _{P_l}[V]^{-2}\Big )(V-\bar{V}),
\end{split}
\end{equation}
along the solutions to \eqref{eq4}, \eqref{eq11}, \eqref{eq16}. Then, we can conclude that $\mathcal{L}S_\mathrm{Z^\ast IP}(x) \leq  (u-{u^\ast})^\top  (I_g-\bar{I}_g)$ for all trajectories $(I_g,V,I,\hat{I}_l,\hat{P}_l)\in\Omega_\mathrm{Z^\ast IP}$.
\hfill $\blacksquare$
\end{pf}

\subsection{ZIP loads}
In this subsection, we finally consider ZIP loads, where the SDE describing the dynamics of $\hat{G}_{li}$ is given by 
\begin{equation} \label{eq20}
d\hat{G}_{li}=-\mu _{G_{li}}\hat{G}_{li} dt + \sigma _{G_{li}} \hat{G}_{li} dw_i,
\end{equation}
where $\mu _{G_{li}}$ and $\sigma _{G_{li}}$ are positive constants. Also in this case, the deterministic and stochastic parts of (\ref{eq20}) at the equilibrium point are identical to zero.
Now, in order to verify a passivity property of the considered network, we introduce an assumption on the parameters of the SDE \eqref{eq20} and restrict the trajectories to be inside a subset of the state-space. 
\begin{assumption}{\bf(Condition on the parameters of \eqref{eq20}).}\label{ass:6}
For all $i\in\mathcal{V}$, the stochastic parameters of (\ref{eq20}) satisfy
\begin{equation}\label{eq17.1}
\mu _{G_{l}} >\sigma _{G_{l}}^2,
\end{equation}
where $\mu _{G_l}=\diag\{\mu _{G_{l1}},\dots,\mu _{G_{ln}} \}$ and $\sigma _{G_l}=\diag\{\sigma _{G_{l1}},$ $\dots,\sigma _{G_{ln}}\}$. 
\end{assumption}

Let us define the set $\Omega_\mathrm{ZIP}:=\{ (I_g,V,I,\hat{I}_l,\hat{P}_l,$ $\hat{G}_l): \frac{1}{2}([G_l^\ast]-\Pi^{-1})-[P^\ast_l][V]^{-1}[\bar{V}]^{-1}-\frac{1}{2}\Sigma\mu _{P_l}[V]^{-2}- [V-\bar{V}]^2\Lambda\mu_{G_l}> 0\}$.
Then, the stochastic passivity property of system (\ref{eq4}) with ZIP loads is verified via the following lemma.  
\begin{lemma}({\bf Stochastic passivity of \eqref{eq4}, \eqref{eq11}, \eqref{eq16}, \eqref{eq20}).}
Let Assumptions \ref{ass:2}--\ref{ass:6} hold and $\Omega_\mathrm{ZIP}$ be nonempty. System (\ref{eq4}), with $\hat{I}_{l}$, $\hat{P}_{l}$, $\hat{G}_{l}$ given by (\ref{eq11}), (\ref{eq16}), (\ref{eq20}) is stochastically (shifted) passive with respect to the storage function
\begin{equation}\label{eq18}
\begin{split}
S_\mathrm{ZIP}(I_g, V, I, \hat{I}_l,\hat{P}_l,\hat{G}_l)=&~S_\mathrm{Z^\ast IP}(I_g, V, I, \hat{I}_l,\hat{P}_l)+\frac{1}{2}\hat{G}_l^\top\Lambda \hat{G}_l,
\end{split}
\end{equation}
and supply rate $(I_g-\bar{I}_g)^\top (u-{u^\ast})$ for all the trajectories $(I_g,V,I,\hat{I}_l,\hat{P}_l,\hat{G}_l)\in\Omega_\mathrm{ZIP}$, where $\Lambda$ is a (suitable) positive-definite diagonal constant matrix of appropriate dimensions.
\end{lemma}
\begin{pf}
The Ito derivative of the storage function (\ref{eq18}) satisfies
\begin{equation}\label{eq19.1}
\begin{split}
\mathcal{L}S_\mathrm{ZIP}=&(I_g-\bar{I}_g)^\top (u-{u^\ast})-(I-\bar{I})^\top R(I-\bar{I})\\
&-\hat{I}_l^\top (\mu_{I_l}-\frac{1}{2} \sigma _{I_l}^2-\frac{1}{2}\mathbb{I}_{n})\Pi\hat{I}_l\\
&-\frac{1}{2}(\Pi^{-\frac{1}{2}}(V-\bar{V})+\Pi^{\frac{1}{2}}\hat{I}_l)^\top (\Pi^{-\frac{1}{2}}(V-\bar{V})\\
&+\Pi^{\frac{1}{2}}\hat{I}_l)-\frac{1}{2}\hat{P}_l^\top (\mu _{P_l}-\sigma _{P_l}^2)\Sigma \hat{P}_l\\
&-\frac{1}{2}\Big (\Sigma ^{\frac{1}{2}}\mu _{P_l}^\frac{1}{2}\hat{P}_l-\Sigma ^{-\frac{1}{2}}\mu _{P_l}^{-\frac{1}{2}}(-\mathds{1}_{n}^\top+\bar{V}^\top [V]^{-1})^\top\Big )^\top \\
&\Big (\Sigma ^{\frac{1}{2}}\mu _{P_l}^\frac{1}{2}\hat{P}_l-\Sigma ^{-\frac{1}{2}}\mu _{P_l}^{-\frac{1}{2}}(-\mathds{1}_{n}^\top+\bar{V}^\top [V]^{-1})^\top\Big)\\
&-\frac{1}{2}\hat{G}_l^\top (\mu_{G_l}-\sigma _{G_l}^2)\Lambda \hat{G}_l -\frac{1}{2}\Big(\Lambda ^{\frac{1}{2}}\mu _{G_l}^{\frac{1}{2}}\hat{G}_l\\
&+\Lambda ^{\frac{-1}{2}}\mu _{G_l}^{\frac{-1}{2}}(V-\bar{V})(V-\bar{V})^\top \mathds{1}_{n}\Big )^\top \Big(\Lambda ^{\frac{1}{2}}\mu _{G_l}^{\frac{1}{2}}\hat{G}_l\\
&+\Lambda ^{\frac{-1}{2}}\mu _{G_l}^{\frac{-1}{2}}(V-\bar{V})(V-\bar{V})^\top\mathds{1}_{n}\Big )-(V-\bar{V})^\top \\
&\Big (\frac{1}{2}([G^\ast _l]-\Pi ^{-1})-[P^\ast_l][V]^{-1}[\bar{V}]^{-1}-\frac{1}{2}\Sigma\mu _{P_l}[V]^{-2}\\
&-[V-\bar{V}]^2\Lambda\mu_{G_l}\Big )(V-\bar{V}),
\end{split} 
\end{equation}
along the solutions to \eqref{eq4}, \eqref{eq11}, \eqref{eq16}, \eqref{eq20}. Then, we can conclude that $\mathcal{L}S_\mathrm{Z^\ast IP^\ast}(x) \leq  (u-\bar{u})^\top (I_g-\bar{I}_g)$ for all trajectories $(I_g,V,I,\hat{I}_l,\hat{P}_l,\hat{G}_l)\in\Omega_\mathrm{ZIP}$.
\hfill $\blacksquare$
\end{pf}
\begin{remark}{\bf (Sufficient conditions)}
We observe that for large $\Pi$ and small $\Sigma, \Lambda$, the sufficient conditions for $\Omega_\mathrm{ZIP}$ to be nonempty are similar to the (well-known) sufficient conditions provided in the literature for DC networks with (non-stochastic) Z$^\ast$I$^\ast$P$^\ast$ loads, \textit{i.e.}, high voltage and large values of the load conductance (\cite{ref20.5,ref23,ref22,CucuzzellaCDC19}).
\end{remark}

\subsection{Closed-loop analysis}
In this subsection, we consider the distributed controllers proposed by \cite{ref24.5} and show that the closed-loop system is asymptotically stochastically stable, achieving at the steady state Goals 1 and 2. Before introducing the controller proposed by \cite{ref24.5}, we recall for the readers' convenience the definition of (asymptotic) stochastic stability 
(\cite{ref25,ref26,ref26.1,ref30}).{
\begin{definition} {\bf((Asymptotic) stochastic stability).}  System (\ref{eq7}) is (asymptotically) stochastically stable if a twice continuously differentiable positive definite Lyapunov function $S: \R^N\longrightarrow \R_{>0}$ exists such that $\mathcal{L}S$ is (negative definite) negative semi-definite.
\end{definition}}
Now, for the sake of completeness, we report below the controller proposed by \cite{ref24.5} for $i \in \mathcal{V}$, \textit{i.e.},
\begin{align}\label{eq20.1}
    \tau_{\xi i} \dot \xi _i =& -\sum_{j \in \mathcal{N}^{com}_i}\rho _{ij}(q_i I_{gi} - q_j I_{gj})  \nonumber \\
  \tau_{\eta i} \dot \eta _i =& -\eta_i + I_{gi} \\
  u_i =&~ -K_i(I_{gi} - \eta _i) +q_i \sum_{j \in \mathcal{N}^{com}_i}\rho _{ij}(\xi _i - \xi_j) + V_i^{\ast}. \nonumber
\end{align}
where $\tau_{\xi i}, \tau_{\eta i}, K _i\in \mathbb{R}_{> 0}$ are design parameters and $\mathcal{N}^{com}_i$ is the set of DGUs communicating with DGU $i$ via a communication network (possibly different from the electric network), where $\rho_{ij}\in\R_{>0}$ are the edge weights. 
The distributed controller (\ref{eq20.1}) can be written compactly for all $i\in\mathcal{V}$ as 
\begin{subequations}\label{eq24}
\begin{align}
\tau_{\xi} \dot \xi =& -\mathcal{L}^{com}Q I_g \label{eq24a} \\
  \tau_{\eta} \dot \eta =& -\eta + I_g \label{eq24b} \\
  u =&~ -K(I_g - \eta) +Q \mathcal{L}^{com}\xi + V^{\ast}, \label{eq24c}
\end{align}
\end{subequations}
where $\tau_{\xi}, K, \tau_{\eta}$ are positive definite diagonal matrices of appropriate dimensions and $\mathcal{L}^{com}$ is the weighted Laplacian matrix associated with the communication network. 

In the following theorem, we show that the closed-loop system is asymptotically stochastically stable and the average voltage regulation and current sharing goals are attained.
\begin{theorem}{\bf(Closed-loop analysis).}
Let Assumptions \ref{ass:1}--\ref{ass:6} hold and $\Omega_\mathrm{ZIP}$ be non-empty. System (\ref{eq4}) with $\hat{I}_{l}$, $\hat{P}_{l}$, $\hat{G}_{l}$  given by (\ref{eq11}), (\ref{eq16}), (\ref{eq20}) and controlled by (\ref{eq24}) is asymptotically stochastically stable and achieves Goals 1 and 2 for all the trajectories $(I_g,V,I,\hat{I}_l,\hat{P}_l,\hat{G}_l)\in\Omega_\mathrm{ZIP}$. 
\end{theorem}

\begin{pf}
Let $x:=(I_g,V,I,\hat{I}_l,\hat{P}_l,\hat{G}_l,\xi,\eta)$ and consider the following storage function
\begin{align}\label{eq26}
\begin{split}
S(x)=&~S_\mathrm{ZIP}(I_g, V, I, \hat{I}_l,\hat{P}_l,\hat{G}_l)\\
&+\frac{1}{2}(\xi - \bar{\xi})^\top\tau_{\xi}(\xi -\bar{\xi})+\frac{1}{2}(\eta -\bar{\eta})^\top\tau_{\eta}(\eta -\bar{\eta}),
\end{split}
\end{align}
The Ito derivative of the storage function (\ref{eq26}) satisfies
\begin{equation}\label{eq27}
\begin{split}
\mathcal{L}S=&(I_g-\eta)^\top K(I_g-\eta)-(I-\bar{I})^\top R(I-\bar{I}) \\
&-\hat{I}_l^\top (\mu_{I_l}-\frac{1}{2} \sigma _{I_l}^2-\frac{1}{2}\mathbb{I}_{n})\Pi\hat{I}_l\\
&-\frac{1}{2}(\Pi^{-\frac{1}{2}}(V-\bar{V})+\Pi^{\frac{1}{2}}\hat{I}_l)^\top (\Pi^{-\frac{1}{2}}(V-\bar{V})\\
&+\Pi^{\frac{1}{2}}\hat{I}_l)-\frac{1}{2}\hat{P}_l^\top (\mu _{P_l}-\sigma _{P_l}^2)\Sigma \hat{P}_l\\
&-\frac{1}{2}\Big (\Sigma ^{\frac{1}{2}}\mu _{P_l}^\frac{1}{2}\hat{P}_l-\Sigma ^{-\frac{1}{2}}\mu _{P_l}^{-\frac{1}{2}}(-\mathds{1}_{n}^\top +\bar{V}^\top [V]^{-1})^\top\Big )^\top \\
&\Big (\Sigma ^{\frac{1}{2}}\mu _{P_l}^\frac{1}{2}\hat{P}_l-\Sigma ^{-\frac{1}{2}}\mu _{P_l}^{-\frac{1}{2}}(-\mathds{1}_{n}^\top +\bar{V}^\top [V]^{-1})^\top\Big )\\
&-\frac{1}{2}\hat{G}_l^\top (\mu_{G_l}-\sigma _{G_l}^2)\Lambda \hat{G}_l  -\frac{1}{2}\Big(\Lambda ^{\frac{1}{2}}\mu _{G_l}^{\frac{1}{2}}\hat{G}_l\\
&+\Lambda ^{-\frac{1}{2}}\mu _{G_l}^{-\frac{1}{2}}(V-\bar{V})(V-\bar{V})^\top\mathds{1}_{n}\Big)^\top \Big(\Lambda ^{\frac{1}{2}}\mu _{G_l}^{\frac{1}{2}}\hat{G}_l\\
&+\Lambda ^{-\frac{1}{2}}\mu _{G_l}^{-\frac{1}{2}}(V-\bar{V})(V-\bar{V})^\top\mathds{1}_{n}\Big)-(V-\bar{V})^\top\\
&\Big ( \frac{1}{2}([G^\ast _l]-\Pi ^{-1})-[P^\ast_l][V]^{-1}[\bar{V}]^{-1}-\frac{1}{2}\Sigma\mu _{P_l}[V]^{-2}\\
&-[V-\bar{V}]^2\Lambda\mu_{G_l}\Big )(V-\bar{V}),
\end{split}
\end{equation}
along the solutions to the closed-loop system \eqref{eq4}, \eqref{eq11}, \eqref{eq16}, \eqref{eq20}, \eqref{eq24}. Then, it follows that $\mathcal{L}S \leq 0$, for all the trajectories $(I_g,V,I,\hat{I}_l,\hat{P}_l,\hat{G}_l)\in\Omega_\mathrm{ZIP}$. 
Then, as a preliminary result we can conclude that the solutions to the closed-loop system (\ref{eq4}), (\ref{eq11}), (\ref{eq16}), (\ref{eq20}), (\ref{eq24}) are bounded.
Moreover, according to LaSalle's invariance principle, these solutions converge to the largest invariant set contained in $\Psi:=\{I_g,I,V,\hat{I}_l,\hat{P}_l,\hat{G}_l,\xi ,\eta : I_g=\eta , I=\bar{I}, V=\bar{V}, \hat{I}_l=\hat{P}_l=\hat{G}_l=\textbf{0}\}$. Hence, the behavior of the closed-loop system (\ref{eq4}), (\ref{eq11}), (\ref{eq16}), (\ref{eq20}), (\ref{eq24}) on the set $\Psi$ can be described by
\begin{subequations}
\begin{align}
L_g\dot{I}_{g}  =&- \bar{V}  +Q \mathcal{L}^{com}\xi + V^{\ast} \label{seta}\\
\textbf{0}  = &~I_{g} + \mathcal{A}\overline I - I^\ast_l -[G^\ast _l]\bar{V}-[\bar{V}]^{-1}P^\ast_l\\
\boldsymbol{0} = &-\mathcal{A}^\top \bar{V} - R \overline I \label{setc}\\
d\hat{I}_l  =&~\textbf{0}\\
d\hat{P}_l =&~\textbf{0}\\
d\hat{G}_l =&~\textbf{0}\\
\tau_{\xi} \dot \xi =& -\mathcal{L}^{com}Q I_g \label{setd}\\
\tau_{\eta }\dot \eta =& ~\boldsymbol{0} \label{sete}.
\end{align}
\end{subequations}

It follows from (\ref{sete}) that $\eta$ is constant on the largest invariant set. Then, $I_g$ is also constant. Since $I_g$ is constant, it follows from (\ref{setd}) that $\mathcal{L}^{com}Q I_g$ is also constant, implying that $\xi$ would increase unbounded if $\mathcal{L}^{com}Q I_g\neq\textbf{0}$, contradicting the preliminary result on the boundedness of the solutions.
Then, $\mathcal{L}^{com}Q \bar{I}_g$ must necessarily be equal to zero, implying that the goal of current sharing is achieved (see Goal 1).
Since $I_g$ is constant, we can pre-multiply (\ref{seta}) by $\mathds{1}_{n}^\top Q^{-1}$ and obtain $\mathds{1}_{n}^\top Q^{-1} \bar{V} = \mathds{1}_{n}^\top Q^{-1} V^{\ast}$, \textit{i.e.},  average voltage regulation (see Goal 2). 
\hfill $\blacksquare$
\end{pf}

Note that, in analogy with Theorem~1 and by virtue of Lemmas 1 and 2, similar results can be proved for the cases of Z$^\ast$IP$^\ast$ and Z$^\ast$IP loads, respectively.
\section{Simulation Results}
\begin{figure}[t]
\centering
\includegraphics[width=5cm, height=4cm]{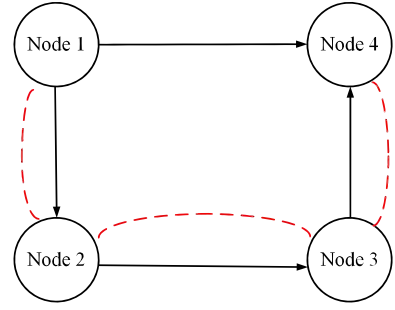}
\caption{Topology of DC microgrid. The arrows denote the positive direction of $I_g$; the dashed lines represent the communication links.}
\label{f2}
\end{figure}
\begin{figure}[t]
\centering
\includegraphics[width=\columnwidth]{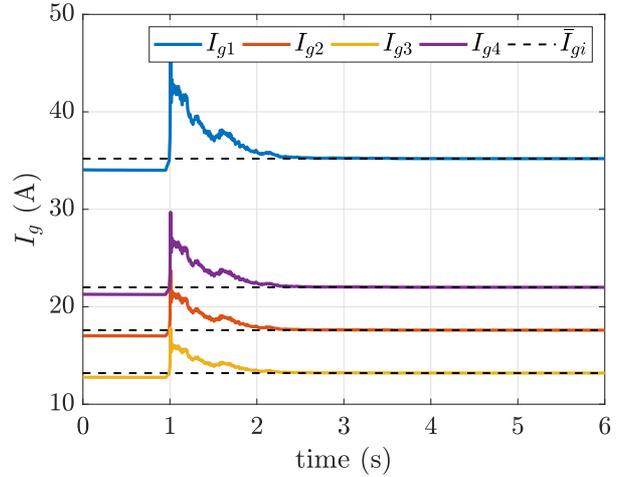}
\caption{The generated current at different nodes.}
\label{f3}
\end{figure}
\begin{figure}[t]
\centering
\includegraphics[width=\columnwidth]{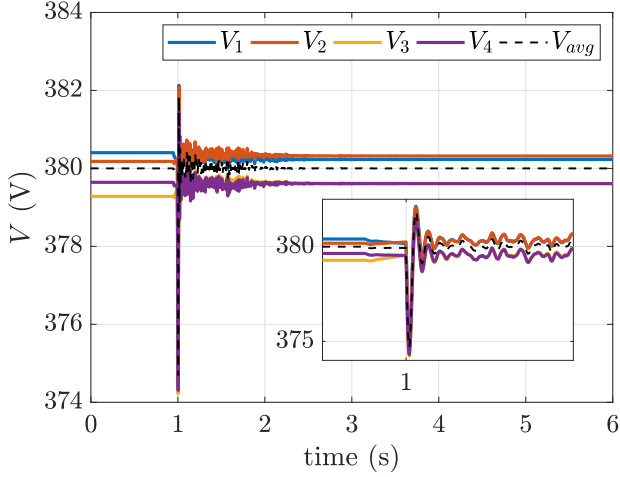}
\caption{Voltage of PCC at different nodes.}
\label{f4}
\end{figure}
\begin{figure}[t]
\centering
\includegraphics[width=\columnwidth]{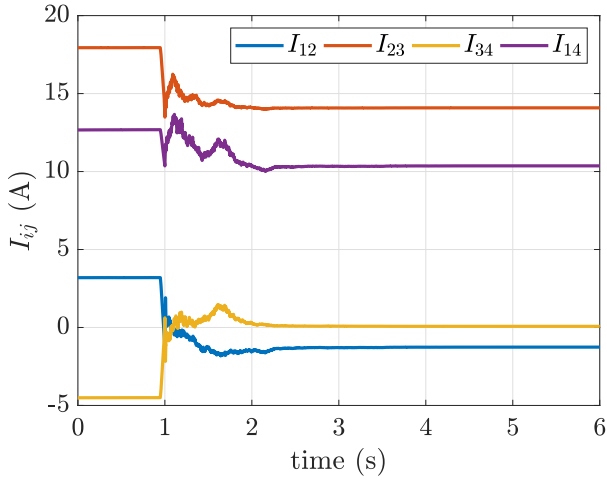}
\caption{The line currents.}
\label{f5}
\end{figure}
\begin{figure}[t]
\centering
\includegraphics[width=\columnwidth]{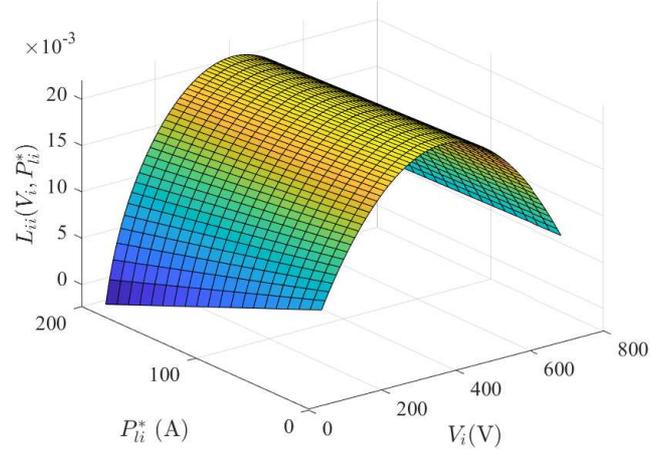}
\caption{The simulation of $L_{ii}(V_i,P^\ast_{li})$ for different $V_i$ and $P^\ast_{li}$.}
\label{f6}
\end{figure}
In this section, the performance of the proposed method is evaluated via a case study example. We consider the closed-loop system (\ref{eq4}), (\ref{eq11}), (\ref{eq16}), (\ref{eq20}), (\ref{eq24}) describing a DC network composed of 4 nodes (\textit{i.e.}, DGUs) depicted in Fig.~\ref{f2}. The parameters of the DC network are taken from \cite{ref12}. The parameters of the controllers are chosen as $K_i=0.4, \tau_{\eta i}=0.005, \tau_{\xi i}=1, i=1,..., 4$. The initial values of $P^\ast_l$ are $\left[25, 10, 25, 20\right]$. Their variations at the time instant $t=1s$  are $\Delta P^\ast_l=\left[8, 5, -8, 3\right]$, while the initial values of $I^\ast_l$, and $G^\ast_l$ are $\left[0.07, 0.045, 0.06, 0.08\right]$ and $\left[8, 4, 5, 12\right]$. The stochastic parameters related to the load current, power and conductance are selected as $\mu_{I_{li}}=2.5, \sigma_{I_{li}}=1, \mu_{P_{li}}=2, \sigma_{P_{li}}=0.7, \mu_{G_{li}}=1.3, \sigma_{G_{li}}=0.2, i=1,..., 4$, respectively. Fig.~\ref{f3} illustrates the  currents generated by each DGU. We can see that the generated currents converge to the values that ensure fair current sharing (dashed lines). Moreover, Fig.~\ref{f4} shows the voltage at the PCC of each node. It can be seen that the weighted average voltage (dashed line) converges to the weighted average of the desired voltages, \textit{i.e.}, 380 V. The current flows through the transmission lines are represented in Fig.~\ref{f5}.

{To investigate the size in practice of the set $\Omega_\mathrm{ZIP}$, we consider $L(V,P^\ast_l):=\frac{1}{2}([G_l^\ast]-\Pi^{-1})-[P^\ast_l][V]^{-1}[\bar{V}]^{-1}-\frac{1}{2}\Sigma\mu _{P_l}[V]^{-2}- [V-\bar{V}]^2\Lambda\mu_{G_l}$, with $\Pi=10^3 \mathbb{I}_{n}$, $\Sigma=\Lambda=10^{-7} \mathbb{I}_{n}$, $[G_l^\ast]=0.045 \mathbb{I}_{n}$. Fig.~\ref{f6} illustrates $L_{ii}(V_i,P^\ast_{li})$  for different values of $V_i$ and $P^\ast_{li}$, where $i=2$. It can be seen that $L_{ii}(V_i,P^\ast_{li})>0$ for a wide range of $V_i$ and $P^\ast_{li}$, \textit{i.e.}, $V_i\in[60,800]$, and $P^\ast_{li}\in[5,200]$. Therefore, the stability of the system is guaranteed for a wide range of operating conditions. Moreover, we notice that by increasing the value of $\Pi_{ii}$ or decreasing the values of $\Sigma_{ii}$ and $\Lambda_{ii}$, the range of $V_i$ and $P^\ast_{li}$ can be enlarged. Furthermore, we notice that these constant matrices are parameters of the Lyapunov function and do not affect the parameters of the used controller.}

\section{CONCLUSIONS}
{In this paper, we have considered stochastic dynamics for the impedance (Z), current (I) and power (P) components of ZIP loads in DC power network. Then, we have verified the stochastic passivity of the considered system. In order to achieve {average} voltage regulation and current sharing, we have used an existing distributed control scheme, proving the asymptotic stochastic stability of overall system. Future research includes the modeling of the ZIP loads as exo-systems and using output regulation methods to tackle the problem of {(average)} voltage regulation and current sharing in DC networks.}

\balance

\end{document}